\documentstyle[12pt]{amsart}
\theoremstyle{plain}
\newtheorem{theorem}{Theorem}
\newtheorem{lemma}{Lemma}
\newtheorem{proposition}{Proposition}
\newtheorem{corollary}[proposition]{Corollary}

\theoremstyle{definition}
\newtheorem*{definition}{Definition}
\newtheorem*{notation}{Notation}

\theoremstyle{remark}
\newtheorem{remark}{Remark}
\newtheorem*{assumption}{Assumption}

\numberwithin{equation}{section}

\newcommand{\ep}{\varepsilon}
\newcommand{\ffi}{\varphi}

\newcommand{\m}{\text{m}}
\newcommand{\C}{\mathbb C}
\newcommand{\I}{\mathbb I}

\newcommand{\R}{\mathbb R}
\newcommand{\T}{\mathbb T}
\newcommand{\Z}{\mathbb Z}
\newcommand{\U}{\text{U}}
\newcommand{\cA}{\mathcal A}
\newcommand{\cB}{\mathcal B}
\newcommand{\cD}{\mathcal D}
\newcommand{\cF}{\mathcal F}
\newcommand{\cH}{\mathcal H}

\newcommand{\coker}{\text{coker}\,}

\newcommand{\grad}{\text{grad}\;}
\newcommand{\ima}{\text{im}\,}

\newcommand{\dx}{\frac{\partial}{\partial x}}
\newcommand{\dy}{\frac{\partial}{\partial y}}
\newcommand{\dt}{\frac{\partial}{\partial\theta}}
\newcommand{\df}{\frac{\partial}{\partial\ffi}}
\begin{document}

\title{Reidemeister torsion and Integrable Hamiltonian systems}
\author[A. Fel'shtyn]{Alexander Fel'shtyn}  
\author[H. S\'anchez-Morgado]{Hector S\'anchez-Morgado}

\maketitle


\setcounter{section}{-1}
  \section{Introduction}

Reidemeister torsion is a very important topological invariant which 
has useful applications in knot theory, quantum field theory and dynamical 
systems. In 1935 Reidemeister \cite{re} classified up to $PL$ equivalence the 
lens spaces $S^3/\Gamma$ where $\Gamma$ is a finite cyclic group of fixed 
point free orthogonal transformations. He used a certain new invariant 
- the Reidemeister torsion- which was quickly extended by Franz, 
who used it to classify the generalized lens spaces $S^{2n+1}/\Gamma$.
Let $X$ be a compact smooth manifold. A representation 
$\rho:\pi_1(X)\rightarrow\U(\m)$ of the fundamental group defines a flat 
${\C}^m$ bundle $E$ over $X$. When the twisted cohomology $H^*(X;E)$ vanishes, 
the representation $\rho$ and the flat bundle $E$ are called acyclic.
The Reidemeister torsion or $R$-torsion is a positive number which
is a ratio of determinants concocted from the $\pi_1(X)$-equivariant
chain complex of the universal covering of $X$.  Later Milnor identified the 
Reidemeister torsion with the Alexander polynomial, 
which plays a fundamental role in the theory of knots and links.

In 1971, Ray and Singer \cite{rs} introduced an analytic torsion  associated 
with the de Rham complex of forms with coefficients in a flat bundle over a 
compact Riemannian manifold, and conjectured it was the same as the 
Reidemeister torsion. The Ray- Singer conjecture was established independently
by Cheeger and M\"uller a few years later.

Recently, the Reidemeister torsion has found interesting applications in 
dynamical systems theory. A connection between the Lefschetz type dynamical 
zeta functions and the Reidemeister torsion was established by D. 
Fried \cite{fri}. The work of Milnor \cite{mi} was the first indication that 
such a connection exists.

In this paper we study the Reidemeister torsion of isoenergy surfaces of an
integrable Hamiltonian system. 
Let $N$ be a four-dimensional smooth symplectic manifold and consider 
the Hamiltonian system with smooth Hamiltonian $H$, which
in Darboux coordinates has the form:

\begin{equation}\label{hamilton}
  \begin{split}
\frac{dp_i}{dt} & = -\frac{\partial H}{\partial q_i}\\
\frac{dq_i}{dt} & = \frac{\partial H}{\partial p_i}.
   \end{split}
\end{equation}

The three-dimensional level surface $M=\{H=const\}$  is invariant under the 
flow defined by the system \eqref{hamilton}. The surface $M$ is called an
isoenergey surface or a constant-energy surface. The topological structure 
of isoenergy surfaces of integrable Hamiltonian systems and the structure of 
their fundamental groups were describeed in \cite{fo, Ziesch}. Isoenergy 
surfaces of integrable Hamiltonians system possess specific properties 
which distinguish them among all smooth three dimensional manifolds.
Namely, they belong to the class of graph-manifolds introduced by 
Waldhausen\cite{wal}. Since $N$ is orientable (as a symplectic manifold), 
the surface $M$ is automatically orientable in all cases. Suppose that the 
system \eqref{hamilton} is complete integrable (in Liouville's sense) on the 
surface $M$. This means, that there is a smooth function $f$
(the second integral), which is independent of $H$ and
with Poisson bracket $\{H,f\}=0 $ in a neighborhood of $M$. 
\begin{definition}
We shall call $f:M\to\R$ a Bott function if its critical points form 
critical nondegenerate smooth submanifolds of $M$. This means that the Hessian 
$d^2f$ of the function $f$ is nondegenerate on the planes normal to the 
critical submanifolds of the function $f$. 
\end{definition}

A.T. Fomenko \cite{fo} proved that a Bott integral on a compact nonsingular 
isoenergy surface $M$ can have only three types of critical submanifolds: 
circles, tori or Klein bottles. The investigation of concrete mechanical 
and physical systems shows \cite{fo} that it is a typical situation when the 
integral on $M$ is a Bott integral. In the classical integrable cases of the 
solid body motion (cases of the Kovalevskaya, Goryachev-Chaplygin, Clebsch,
Manakov) the Bott integrals are round Morse functions on the isoenergy 
surfaces. A round Morse function is a Bott function all whose critical
manifolds are circles. Note that critical circles of $f$ are periodic 
solutions of the system \eqref{hamilton} and the number of this circles is 
finite. Suppose for the moment that the Bott integral $f$ is a round Morse 
function on the closed isoenergy surface $M$. Let us recall the concept of 
the  separatrix diagram of the critical circle $\gamma$. 
Let $x\in \gamma $ be an arbitrary point and  $N_x(\gamma)$ be a disc of 
small radius normal to $\gamma$ at $x$. The restriction of $f$ to the 
$N_x(\gamma)$ is a normal Morse function with the critical point $x$ having a 
certain index $u(\gamma)=0,1,2$. A separatrix of the critical point $x$ is an
integral trajectory of the field $-\grad f$, called a gradient line, which is 
entering  or leaving $x$. The union of all the separatrices leaving the point
$x$ gives a disc of dimension $u(\gamma)$ and is called the outgoing 
separatrix diagram (disc). The union of incoming separatrices gives a 
disc of complementary dimension and is called the separatrix incoming diagram 
(disc). Varing the point $x$ and constructing the incoming and outgoing
separatrix discs for each point $x$, we obtain  the incoming and outgoing 
separatrix diagrams of the circle $\gamma$. Let $\Delta(\gamma)$ be +1 if 
the outgoing separatrix is orientable, and -1 if it is not. 
Let $\epsilon(\gamma)=(-1)^{u(\gamma)}$. 
Let $\rho_E: \pi_1(M,p) \rightarrow\U(E_p)$ be the holonomy representation of 
the hermitian bundle $E$ over $M$; $E_p $ is the fiber at the point $p$.
For the gradient flow of $f$ on $M$, one can construct an {\em index
filtration} for the collection of critical circles $\{\gamma_i\}$, i. e.
a collection of compact submanifolds $M_i$ of top dimension so that 
$M_i \subset\text{int}M_{i+1}, M_0=\emptyset, M_i=M$ for large $i$, the flow
is transverse inwards on $\partial M_i$ and $M_{i+1}\setminus M_i$ is an
isolating neighborhood for the critical circle $\gamma_i$.
When  $E|(M_{i+1}, M_i)$ is acyclic for each $i$, following the ideas of 
D. Fried \cite{fri} we can compute the Reidemeister torsion as \cite{fel}
\begin{equation}
\tau(M;E)=\prod_i \tau(M_{i+1}, M_i;E)=
\prod_{\gamma_i}|\det(I-\Delta(\gamma_i)\cdot\rho_E(\gamma_i))|^{\ep(\gamma_i)}
\end{equation}
This formula means that for the integrable Hamiltonian system 
on the four-dimensional symplectic manifold, the Reidemeister torsion of the 
isoenergy surface counts the critical circles  of the second independent 
Bott integral on this surface. If  $E|(M_{i+1}, M_i)$ is acyclic, then 
$\det(I-\Delta(\gamma_i)\cdot\rho_E(\gamma_i))\ne 0$ for each $i$.
Since in many  classical integrable cases there are contractible critical 
circles it is interesting to study the situation when not all 
$E|(M_{i+1}, M_i)$ are acyclic. In this paper we carry out this study and in 
fact we consider the general situation when the Bott integral has critical 
tori and Klein bottles. We use the spectral sequence defined by the filtration
and following Witten-Floer  ideas we bring into play the orbits connecting the 
critical submanifolds. A similar approach was developed in \cite{sm} 
for Morse-Smale flows.

Parts of this article  were written while the first  author was visiting the Instituto de Matematicas Universidad Nacional Autonoma de Mexico  and Erwin Schr\"odinger International Institute for Mathematical Physics in Wien .The first  author is indebted to these institutions for their invitations, support  and 
hospitality.

\section{R-torsion and Spectral Sequences}\label{spectral}
Let $W$ be a finite dimensional vector space with basis
${\bf w}=\{w_1,\ldots,w_n\}$, then $\wedge{\bf w}=w_1\wedge\cdots\wedge w_n$
is a generator of $\det W=\wedge^nW$. 
If dim $V=0$ set $\det V={\C}$. 

Consider a cochain complex of finite dimensional vector spaces

\begin{equation}\label{cochain}
0\to V^0\buildrel d\over\longrightarrow V^1\to\cdots\to V^m
\buildrel d\over\longrightarrow 0 
\end{equation}

Let $V^+=\bigoplus_iV^{2i}$, $V^-=\bigoplus_iV^{2i+1}$ and 

$$\det V=\det(V^-)\otimes(\det V^+)^{-1}.$$

Let $Z^{\pm}=V^{\pm}\cap\ker d$, $B^{\mp}=d(V^{\pm})$,
$H^{\pm}=Z^{\pm}/B^{\pm}$.

We now define the {\em torsion element} 
$\tau_d\in\det V\otimes(\det H)^{-1}$. Pick ordered relative bases
${\bf h}_{\pm}$ for $(Z^{\pm},B^{\pm})$ and 
${\bf t}_{\pm}$ for $(V_\pm,Z_\pm)$,
then $d{\bf t}_{\mp}$ is a basis for $B_{\pm}$. Denote by $[{\bf h}_{\pm}]$ the
corresponding basis for $H^{\pm}$.
\begin{equation}\label{torsion-elem}
\tau_d=\wedge({\bf t}_-,{\bf h}_-,d{\bf t}_+)
\otimes\wedge(d{\bf t}_-,{\bf h}_+,{\bf t}_+)^{-1}\otimes\wedge[{\bf h}_+]
\otimes\wedge[{\bf h}_-]^{-1}
\end{equation}

\begin{notation} Consider the cochain complex
$$0\to V\buildrel A\over\longrightarrow W\to 0$$
We have $H^0=\ker A$, $H^1=\coker A$ and denote $\tau(A):=\tau _d$.
When $A$ is an isomorphism, $\tau(A)$ is the coordinate free version of
$\det A$.
\end{notation}

\begin{proposition} \cite{Freed}
Let $0=F_{N+1}^i\subset F_N^i\subset\cdots\subset F_0^i=V^i$ be a filtration 
of the cochain complex (\ref{cochain}) such that $d^i(F_n^i)\subset 
F_n^{i+1}$. Let $\{E_r,d_r\}$ be the corresponding spectral sequence. Then
$$\tau_d=\tau_{d_0}\otimes\cdots\otimes\tau_{d_N}$$
\end {proposition}

\begin{corollary}\label{aditive}\cite{Freed}
Suppose 
\begin{equation}\label{complex}
0\to (C',d')\buildrel i\over\longrightarrow (C,d)
\buildrel j\over\longrightarrow (C'',d'')\to 0 
\end{equation}
is an exact sequence of chain complexes and
\begin{equation}\label{sequence}
{\cH}: 0\to H^0(C')\buildrel i^*\over\longrightarrow H^0(C)
\longrightarrow H^0(C'')
\buildrel\partial\over\longrightarrow H^1(C')\to\cdots
\end{equation}
is the corresponding long exact sequence. For each $k$ choose compatible
volume elements in $\det C'_k, \det C_k, \det C''_k$, i.e. such that the
torsion of \eqref{complex} is $1$. Then
\begin{equation}
  \label{eq:aditive}
  \tau_d=\tau_{d'}\tau_{d''}\tau_{\cH}
\end{equation}
\end{corollary}
We now describe the first terms of the spectral sequence $\{E_r,d_r\}$.
The filtration defines the associated graded complex $G^i=\bigoplus_nG_n^i$
where $G_n^i=F_n^i/F_{n+1}^i$. The coboundary $d$ induces a map
$d_0^n:G_n^i\to G_n^{i+1}$ whose cohomology defines the term $E_1$ by
\begin{equation}\label{eq:E1}
E_1^{n,q}:=H^{n+q}(F_n/F_{n+1})
\end{equation}
and induces the first differential $d_1^n:E_1^{n,q}\to E_1^{n+1,q}$ as
the coboundary map for the short exact sequence

$$0\to F_{n+1}/F_{n+2}\to F_n/F_{n+2}\to F_n/F_{n+1}\to 0$$
i.e. the map $d_1$ in the long exact sequence

\begin{equation}\label{long-sequence}
\buildrel j_n\over\longrightarrow H^{n+q}(F_n/F_{n+2}) 
\buildrel k_n\over\longrightarrow H^{n+q}(F_n/F_{n+1})
\buildrel d_1\over\longrightarrow H^{n+q+1}(F_{n+1}/F_{n+2})
\end{equation}

The term $E_2^{n,q}$ is defined by
\begin{equation}\label{E2}
E_2^{n,q}:=\frac{\ker(d_1:E_1^{n,q}\to E_1^{n+1,q})}
           {\hbox{im}(d_1:E_1^{n-1,q}\to E_1^{n,q})}.
\end{equation}
From (\ref{long-sequence}) we have $\ker(d_1)=$im$(k_n)$ and 
im$(d_1)=\ker(j_{n-1})$, and thus
$E_2^{n,q}=\hbox{im}(k_n)/\hbox{ker}(j_{n-1}).$

Consider the commutative diagram
$$\begin{array}{rlrlc}
H^{n+q}(F_n/F_{n+2})&\buildrel k_n\over\longrightarrow & H^{n+q}(F_n/F_{n+1})
& \buildrel j_{n-1}\over\longrightarrow & H^{n+q}(F_{n-1}/F_{n+1}) \\
\delta_0  & \nearrow  & \delta_1 & \nearrow & \\
H^{n+q-1}(F_{n-2}/F_n) & \buildrel k_{n-2}\over\longrightarrow &
H^{n+q-1}(F_{n-2}/F_{n-1}) &\buildrel j_{n-3}\over\longrightarrow &
H^{n+q-1}(F_{n-3}/F_{n-1}), 
\end{array}$$
The second differential $d_2^{n-2}:E_2^{n-2,q+1}\to E_2^{n,q}$ is given as
the composite map

\begin{equation}\label{eq:second-diff}
\hbox{im}(k_{n-2})/\hbox{ker}(j_{n-3})
\buildrel\delta_1\over\longrightarrow\hbox{im}(j_{n-1})
\buildrel j_{n-1}^{-1}\over\longrightarrow
\hbox{im}(k_n)/\hbox{ker}(j_{n-1})
\end{equation}

Further terms of the spectral sequence $E_r^{n,q}$ are obtained as
cohomology of the previous term and the differentials
$d_r^n:E_r^{n,q}\to E_r^{n+r,q+r-1}$ are the maps induced by the
original $d$.

Let now $K$ be a finite CW-complex. 
Let $p:\tilde{K}\to K$ be the universal covering and $\rho:\Gamma\to\U(\m)$
be a representation of the fundamental group $\Gamma$  
of $K$ which defines a flat vector bundle 
$E:=\tilde{K}\times_\Gamma{\C}^m$. Lifting
cells to $\tilde{K}$ we obtain a $\Gamma$-invariant CW complex structure
on $\tilde{K}$. The space of $\rho$-equivariant cochains

$$C^\ast(K;E)=\{\xi\in C^\ast(\tilde{K};{\C}^m):
\xi\circ\gamma=\rho(\gamma)\circ\xi\quad\forall\gamma\in \Gamma\}$$
is preserved by $d^j:C^j(\tilde{K};{\C}^m)\to C^{j+1}(\tilde{K};{\C}^m)$
and so $\{C^\ast(K;E), d(K;E)\}$ forms a subcomplex. Its cohomology
$H^\ast(K;E)$ is called the $\rho$-twisted cohomology of $K$.
As usual $H^\ast(K;E)$ is subdivision invariant and we have a torsion 
element

$$\tau_{d(K;E)}\in\det C^\ast(K;E)\otimes(\det H^\ast(K;E))^{-1}.$$

Order the $j$-cells $\sigma$ and choose an oriented lift $\tilde{\sigma}$ for
each $\sigma$. 
This gives an isomorphism  $C^j(K;E)\cong\oplus_{\sigma}{\C}^m$ and 
determines a preferred generator $w_K^\rho$ of $\det(C^\ast(K;E))$ up 
to multiplication by an element of the subgroup

$$U_\rho=\{(\pm 1)^m\det\rho(\gamma):\gamma\in \Gamma\}\subset S^1$$

The orbit $U_\rho w_K^\rho\subset\det C^\ast(K;E)$
is invariant under subdivision, so we can define {\em R-torsion} of $K$ at
$\rho$ as the $U_\rho$ orbit

\begin{equation}\label{r-torsion}
\tau(K;E)=(U_\rho w_K^\rho)^{-1}\otimes\tau_{d(K;E)}\subset
(\det H^\ast(K;E))^{-1}
\end{equation}
which is invariant under subdivision.
When $\rho$ is acyclic, i.e. when $H^\ast(K;E)=0$, we have
$\det H^\ast(K;E)={\C}$ and we can identify $\tau(K;E)$ as an
element of ${\C}^\ast/U_\rho$. Since $U_\rho\subset S^1$,
all elements in $\tau(K;E)$ have the same modulus which we still denote by
$\tau(K;E)$.

The previous definitions can be extended to relative pairs. Let $L$ be
a subcomplex of K. For each $j$ we have the relative space of cellular
$j$-cochains

$$C^j(K,L;{\C})=\bigoplus_{\sigma\in K\setminus L}
H^j(\sigma,\partial\sigma;{\C})$$

Let $\tilde{K}$ and $\rho$  be as above and let $\tilde{L}=p^{-1}(L)$.
We can define the space of relative $\rho$-equivariant cochains
$C^\ast(K,L;E)\subset C^\ast(\tilde{K},\tilde{L};{\C}^m)$ with
coboundary $d(K,L;E)$ and then we get a torsion element

$$\tau_{d(K,L;E)}\in\det C^\ast(K,L;E)\otimes(\det H^\ast(K,L;E))^{-1}.$$

Thus, choosing preferred basis as before we obtain a $U_\rho$ orbit

$$\tau(K,L;E)\subset\det H^\ast(K,L;E)^{-1}$$
which is invariant under subdivision.

\begin{remark} 
Another name for the twisted cohomology is cohomology with 
local coefficients. One chooses a point on each cell of $K$ and a path from
a fixed point to each chosen point. In this way any path $c$ between 
chosen points defines a closed path $\gamma_c$ and then a matrix 
$\rho(c):=\rho(\gamma_c)$ which gives the relation between the 
coefficients at the ends of the path.
\end{remark}

\section{R-torsion and critical submanifolds}
The foliation of the isoenergy surface $M$ by Liouville tori is given
by the level sets of the Bott integral $f^{-1}(c)$ for $c$ a regular
value. The bifurcation of Liouville tori ocurr at the sets 
$F_c=f^{-1}(c)\cap\text{Crit}(f)$ for $c$ a critical value. 
We will make the following assumption wich is satisfied in the generic case.
 \begin{assumption}
We will assume that there are no gradient lines of the Bott integral 
$f$ connecting saddle circles i.e. circles with index $1$.
 \end{assumption}
We will substitute the Bott integral for another Bott function ,
still denoted by $f$ and not necessarily an integral,
giving the same foliation by Liouville tori and such that its critical values
$c_1<\cdots<c_l$, are ordered in the following way 
\begin{itemize}
\item [(a)]$i\le k_1 \iff F_{c_i}$ is a minimum circle.
\item [(b)]$ k_1<i\le k_2 \iff F_{c_i}$ is a minimum 
torus or Klein bottle.
\item [(c)]$ k_3<i\le k_4 \iff F_{c_i}$ is a maximum 
torus or Klein bottle.
\item [(d)]$ k_4<i \iff F_{c_i}$ is a maximum circle.
\end{itemize}
Choosing numbers $A_0<c_1<A_1<\cdots<c_l<A_l$ 
and letting $N_j=f^{-1}(\infty,A_j]$ we get an {\em index filtration}
$\emptyset=N_0\subset N_1\subset\cdots\subset N_l=M$ for the
critical sets of $f$ i. e. denoting by $\phi$ the gradient flow of $f$, 
$\phi$ is transverse inwards on $\partial N_n$ and 
\begin{equation*}
\bigcap_{t\in{\R}}\phi_t(N_j\setminus N_{j-1})=F_{c_j}.
\end{equation*}
Fix a representation $\rho:\pi_1(M)\to\U(\m)$. 
All cochain complexes and cohomology groups will have coefficients in the
flat bundle defined by $\rho$. Let $l_0=0, l_1=k_2, l_2=k_3, l_3=l$ and
$M_n=N_{l_n}$. Define the filtration $F_5\subset\cdots\subset F_1\subset F_0$
of $C^\ast(M)$ by $F_n=$ ker$(C^\ast(M)\to C^\ast(M_n))$.
The associated graded complex is given by
$$G_n=F_n/F_{n+1}\cong C^\ast(M_{n+1},M_n)$$
with 0-differential 
$d_0=\bigoplus_nd_0^n$, where $d_0^n:G_n\to G_n$,
and torsion element $\tau_{d_0}=\bigotimes_n\tau_{d_0^n}$ where
$$\tau_{d_0^n}\in
\det C^\ast(M_{n+1},M_n)\otimes(\det H^\ast(M_{n+1},M_n))^{-1}.$$
Since there are neither gradient lines connecting two minimum (maximum)
critical submanifolds nor gradient lines connecting two saddle circles
(by assumption), we have (see \cite{sm})

\begin{equation}\label{sum}
H^*(M_{n+1},M_n)= \bigoplus_{M_n\subset N_j\subset M_{n+1}}H^*(N_j,N_{j-1})
\end{equation}
\begin{equation}\label{product}
\tau_{d_0^n}=\bigotimes_{M_n\subset N_j\subset M_{n+1}}\tau_{d_0^{n,j}}
\end{equation}

From \eqref{sum}, the computation of the map
$$d_1^n:H^\ast(M_n,M_{n-1})\to H^{\ast+1}(M_{n+1},M_n)$$
reduces to computing for each $i,j$ with $l_{n-1}<i\le l_n<j\le l_{n+1}$, 
its component
$F_{ij}^\ast:H^\ast(N_i,N_{i-1})\to H^{\ast+1}(N_j,N_{j-1})$.
To do so, we will use the trajectories of the gradient flow of $f$,
but we will modify $f$ in the neighborhood of each critical level set in
order to apply Lemma \ref{Witten} below giving such a map in the Morse 
function case. This modification is just a technical device to choose some of
the orbits connecting the critical submanifolds to describe the maps 
$F_{ij}^\ast$. We will give a proof of the following Proposition using Lemma
\ref{Witten}.
\begin{proposition}\label{decomposition}
If $ k_{2n}<j\le k_{2n+1}$, let $\gamma_j=F_{c_j}$ and 
${\cD}_j=I-\Delta(\gamma_j)\rho(\gamma_j)$. Then
%\begin{subequations}\label{torsion-circles}
\begin{equation}
H^k(N_j,N_{j-1})=\begin{cases}\ker{\cD}_j & \hbox{if}\quad k=n,n+1\\
                                        0 & \hbox{in other case}.
                  \end{cases}
\end{equation}
\begin{equation}
\tau_{d_0^{n,j}} = \tau({\cD}_j)^{(-1)^n}
\end{equation}
%\end{subequations}
If $ k_{2n-1}<j\le k_{2n}$ and $\alpha_j,\beta_j$ are generators of the 
fundamental group of $F_{c_j}$, let  
${\cD}_j=\left(\begin{smallmatrix}
I-\rho(\alpha_j)\\I\pm\rho(\beta_j)\end{smallmatrix}\right)$, 
${\cD}_j^*=(I\pm\rho(\beta_j),\rho(\alpha_j)-I)$, where 
the $+$ sign ocurrs precisely when $F_{c_j}$ is a Klein bottle. Then
%\begin{subequations}\label{torsion-tori}
\begin{equation}
 H^k(N_j,N_{j-1})=
          \begin{cases}\ker{\cD}_j=\coker{\cD}_j^* & \hbox{if}\quad |k-n|=1\\
                     \ker{\cD}_j\oplus\ker{\cD}_j  & \hbox{if} \quad k=n\\
                                                 0 & \hbox{in other case}
           \end{cases}
\end{equation}
\begin{equation}
 \tau_{d_0^{n-1,j}} =1 
\end{equation}
%\end{subequations}
\end{proposition}

\begin{corollary}\label{E1}
$$E_1^{0,q}=H^q(M_1)=0\quad\text{for}\quad q\ne 0,1,2.$$
$$E_1^{1,q}=H^{q+1}(M_2,M_1)=0\quad\text{for}\quad q\ne 0,1.$$
$$E_1^{2,q}=H^{q+2}(M,M_2)=0\quad\text{for}\quad q\ne 0,1,2.$$
\end{corollary}

Let $G:M\to{\R}$ be a Morse Smale function and let $c_1<\cdots<c_N$ be 
its critical points. For $A_0<c_1<\cdots<c_N<A_N$ and 
$K_i=G^{-1}(\infty,A_i]$ we get a filtration $K_0\subset\cdots\subset K_N$.
The orientation of  $M$ and $\grad G$ define an orientation of $L_a=G^{-1}(a)$
for each regular value $a$. Giving an orientation to the unstable subspace 
$E^u(x)$ for each critical point $x$ of $G$, and using the orientation of 
$M$ we also get an orientation of $E^s(x)$. Then we have orientations of 
$W^u(x)$ and $W^s(x)$. Let $x,y$ be critical points of $G$ of indices 
$n,n+1$ respectively and let $a$ be a regular value with $G(x)<a<G(y)$.
Then $S^u(y)=W^u(y)\cap L_a$ and $S^s(x)=W^s(x)\cap L_a$ are oriented 
transverse submanifolds of $L_a$ with dimensions $n$ and and $2-n$ 
respectively. Therefore $S^s(x)\cap S^u(y)$ is a finite set. For each  
$q\in S^s(x)\cap S^u(y)$ denote by $I_q$ the intersection number.

The proof given by  Floer in \cite{Floer} for the  untwisted version of 
the following Lemma can be readily adapted. The only new ingredients are the 
matrices $\rho(\alpha)$ for nonclosed paths $\alpha$ used to define cohomology 
with local coefficients.

\begin{lemma}\label{Witten}
Let $c_i,c_{i+1}$ be critical points of $G$ with indices $n,n+1$.
For each $q\in S(i,j)=S^s(c_i)\cap S^u(c_j)$ let 
$\alpha_q(t)=\phi_{\cot(\pi t)}(q):t\in[0,1]$. The coboundary
map $d:H^n(K_i,K_{i-1})\to H^{n+1}(K_{i+1},K_i)$ is given by
\begin{equation}
  \label{Floer}
d=\sum_{q\in S(i,j)}I_q\rho(\alpha_q)
\end{equation}
\end{lemma}

To change the Bott integral in a neighborhood of each critical level set,
we will use the folllowing Propositions.

\begin{proposition}\label{circle}
Let $F:M\to\R$ be a Bott function and let $\gamma$ be a critical 
circle of index $n$. Given a small neighborhood $U$ of $\gamma$ there is 
another Bott function $G$ which agrees with $F$ outside $U$ and has 
nondegenerate critical points $w,z\in\gamma$ of indices $n,n+1$ and no other
critical points in $U$.
\end{proposition}
\begin{pf} Let $F(\gamma)=c$.
If $\Delta(\gamma)=1$, there is a tubular neighborhood $U$ of $\gamma$
with coordinates $(x,y)\in B_\ep(0)$, $\theta\in\times S^1={\R}/{\Z}$ 
such that 
$$F(x,y,\theta)= c + \pm x^2 \pm y^2$$
Let $\rho:{\R}\to[0,1]$ be a smooth function with $\rho(t)=0$ if $t>\ep^2$ and
$\rho(t)=1$ for $t<\ep^2/4$. Define
$$G=\begin{cases}
f+\delta\rho(x^2+y^2)\cos\theta & \text{on $U$}\\
f & \text{outside $U$}
\end{cases}$$
On $U$ we have
\begin{multline*}
\grad G(x,y,\theta)=
2(\pm 1+\delta\rho'(x^2+y^2)\cos\theta)\; x\dx\\
+2(\pm 1+\delta\rho'(x^2+y^2)\cos\theta)\; y\dy-
\delta\rho(x^2+y^2)\sin\theta\dt
\end{multline*}
Thus, if $\delta$ is  sufficiently small, the only critical points of $G$ in 
$U$ are $w=(0,0,\pi)$ and $z=(0,0,0)$ with indices $p$ and $p+1$ respectively.
\end{pf}
\begin{proposition}\label{torus}
Let $F:M\to\R$ be a Bott function and let ${\T}$ be a minimal 
(maximal) torus. Given a small neighborhood $U$ of $\T$ there is another 
Bott function $G$ which 
\begin{itemize}
\item[1.] Agrees with $F$ outside $U$.  
\item[2.] Has nondegenerate critical points $p,q,r,s\in \T$ of indices 
$0,1,1,2$ (or $1,2,2,3$) and no other critical points in $U$.
\item[3.] There are no gradient lines either from a critical point of 
index $1$ to any of $q,r,s$ or from a critical point of index $2$ to $s$
(either from $p$ to a critical point of index $1$ or from $p,q,r,s$ to
a critical point of index $2$).
\end{itemize}
\end{proposition}
\begin{pf}
Consider the case of a mimimal torus. Let $F({\T})=c$.
There is a tubular neighborhood $U$ of $\T$ with coordinates 
$(x,\theta,\ffi)\in(-\ep,\ep)\times S^1\times S^1$ such that 
$F(x,\theta,\ffi)= c+x^2.$

For each critical point $z$ of $F$ of index $2$, 
$W^u(z)\cap\{x\}\times S^1\times S^1$ is a curve. For each critical point $w$ 
of $F$ of index $1$, $W^u(w)\cap\{x\}\times S^1\times S^1$ is a point.
Therefore we can choose the coordinates $\theta,\ffi$ such that 
\begin{itemize}
\item [(a)] $(x,0,0)\notin\underset{u(z)=2}{\bigcup}W^u(z)$
\item [(b)]  $\{x\}\times(S^1-\{\pi\})\times\{0\}
\cap\underset{u(z)=1}{\bigcup}W^u(z)=\emptyset$
\item [(c)]$\{x\}\times\{0\}\times(S^1-\{\pi\})
\cap\underset{i(w)=2}{\bigcup}W^u(w)=\emptyset$
\end{itemize}
Let $\rho$ be as in the proof of Proposition \ref{circle}. Define
$$G=\begin{cases}
F+\rho(x^2)(\delta_1\cos\theta+\delta_2\cos\ffi) & \text{on $U$}\\
F & \text{outside $U$}
\end{cases}$$
On $U$ we have
\begin{equation*}
 \begin{split}
\grad G(x,\theta,\ffi)=&
2(\pm 1+\rho'(x^2)(\delta_1\cos\theta+\delta_2\cos\ffi))\; x\dx\\
- & \rho(x^2)(\delta_1\sin\theta\dt+\delta_2\sin\ffi\df)
 \end{split}
\end{equation*}
Thus, if $\delta_1, \delta_2$ are sufficiently small, the only critical points
of $G$ in $U$ are $p=(0,\pi,\pi), q=(0,\pi,0), r=(0,0,\pi), s=(0,0,0)$ 
which have the required indices and
\begin{itemize}
\item [(a)]$W^s(s)\cap U=(-\ep,\ep)\times\{(0,0)\}$
\item [(b)] $W^s(q)\cap U=(-\ep,\ep)\times(S^1-\{\pi\})\times\{0\}$
\item [(c)]$W^s(r)\cap U=(-\ep,\ep)\times\{0\}\times(S^1-\{\pi\})$
\end{itemize}
The proof of the maximal torus case is similar.
\end{pf}
Applying Proposition \ref{circle} to $f$ we obtain a Bott function 
without critical circles, then using Proposition \ref{torus} 
we obtain a Bott function without critical circles and minimal torus 
(Klein bottles) and with no  gradient lines from a critical point of 
index $i$ to a critical point of index $j\ge i$. 
Using Proposition \ref{torus} again we finally get a Morse-Smale function 
$g$ that agrees with $f$ outside a neighborhood of each critical set.
Denote by $C_j(g)$ the set of critical points of $g$ with index $j$. We have
$$C_0(g)=\{w_1,\ldots,w_{k_1},p_{k_1+1},\ldots,p_{k_2}\}$$
$$C_1(g)=\{z_1,\ldots,z_{k_1},q_{k_1+1},r_{k_1+1},\ldots,q_{k_2},r_{k_2},
w_{k_2+1},\ldots,w_{k_3},p_{k_3+1},\ldots,p_{k_4}\}$$
$$C_2(g)=\{s_{k_1+1},\ldots,s_{k_2},z_{k_2+1},\ldots,z_{k_3},q_{k_3+1},
r_{k_3+1},\ldots,q_{k_4},r_{k_4},w_{k_4+1},\ldots,w_l\}$$
$$C_3(g)=\{s_{k_3+1},\ldots,s_{k_4},z_{k_4+1},\ldots,z_l\}$$

\begin{pf} (Proposition \ref{decomposition})
For $k_{2n}<i\le k_{2n+1}$ let $g(w_i)<d_i<g(z_i)$ and 
$N^*_i=g^{-1}(\infty,d_i]$.  Then $H^n(N^*_i,N_{i-1})={\C}^mw_i$, 
$H^{n+1}(N_i,N^*_i)={\C}^mz_i$, and according to Milnor \cite{Whitehead}
$$\tau_{d(N^*_i,N_{i-1})}=1,\tau_{d(N_i,N^*_i)}=1.$$
Using Lemma \ref{Witten} we have the sequence
\begin{equation}\label{circles}
0\to H^n(N_i,N_{i-1})\buildrel{\cA}_i\over\longrightarrow {\C}^m
\buildrel{\cD}_i\over\longrightarrow{\C}^m\buildrel{\cB}_i
\over\longrightarrow H^{n+1}(N_i,N_{i-1})\to 0
\end{equation}
with ${\cD}_i$ defined in the the statement of Proposition \ref{decomposition}. 
Therefore
$$H^n(N_i,N_{i-1})=\ker{\cD}_i,\quad H^{n+1}(N_i,N_{i-1})=\coker{\cD}_i.$$ 
Since sequence \eqref{circles} has torsion
$\tau({\cD}_i)^{(-1)^n}$, we obtain by Corollary \ref{aditive} 
$$\tau_0^{n,i}=\tau_{d(N_i,N^*_i)}\tau_{d(N^*_i,N_{i-1})}\tau({\cD}_i)^{(-1)^n}
=\tau({\cD}_i)^{(-1)^n}.$$

For $k_{2n-1}<i\le k_{2n}$, let
$g(p_i)<c_i'<g(q_i), g(r_i)<c_i''<g(s_i)$ and 
$N_i'=g^{-1}(\infty,c_i'], N_i''=g^{-1}(\infty,c_i'']$. Then
$H^{n-1}(N_i',N_{i-1})={\C}^mp_i$, 
$H^{n-1}(N_i'',N_i')={\C}^mq_i\oplus{\C}^mr_i$, $H^{n+1}(N_i,N_i'')={\C}^ms_i$.
Again from \cite{Whitehead}
$$\tau_{d(N_i',N_{i-1})}=1, \tau_{d(N_i'',N_i')}=1, \tau_{d(N_i,N_i'')}=1.$$

By Lemma \ref{Witten} we have the sequence of the triad 
$(N_i'',N_i',N_{i-1})$:
\begin{equation}\label{minima}
0\to H^{n-1}(N_i'',N_{i-1})\buildrel{\cA}_i\over\longrightarrow{\C}^m
\buildrel{\cD}_i\over\longrightarrow{\C}^m\oplus{\C}^m
\buildrel {\cB}_i\over\longrightarrow H^n(N_i'',N_{i-1})\to 0,
\end{equation}
and the exact sequence of the triad $(N_i,N_i'',N_i')$:
\begin{equation}\label{maxima}
0\to H^1(N_i,N_i')\buildrel{\cA}_i\over\longrightarrow{\C}^m\oplus{\C}^m
\buildrel{\cD}_i^*\over\longrightarrow{\C}^m
\buildrel{\cB}_i\over\longrightarrow H^2(N_i,N_i')\to 0,
\end{equation}
with ${\cD}_i$ and ${\cD}^*_i$ defined in the statement of Proposition 
\ref{decomposition}.
The first part of the sequence of the triad 
$(N_i,N_i'',N_{i-1})$:
$$0\to H^{n-1}(N_i,N_{i-1})\longrightarrow H^{n-1}(N_i'',N_{i-1})
\longrightarrow 0,$$
and sequence \eqref{minima} give
%\begin{subequations}
  \label{bottom}
\begin{equation}
H^{n-1}(N_i,N_{i-1})= H^{n-1}(N_i'',N_{i-1})=\ker{\cD}_i
\end{equation}
\begin{equation}
H^n(N_i'',N_{i-1})=\coker{\cD}_i.
\end{equation}
%\end{subequations}
Using \eqref{minima} and \eqref{maxima}, we have the commutative diagram
\begin{equation*}
  \begin{array}{lccll}
 &{\cD}_i\downarrow\quad &       &       &   \\
 &{\C}^m\oplus{\C}^m   &              &        &             \\
 &{\cB}_i\downarrow\quad & \searrow {\cD}_i^* \quad &     &    \\
0\to H^n(N_i,N_{i-1})\to & H^n(N_i'',N_{i-1})& \buildrel\Delta_i\over
\longrightarrow{\C}^m & \to & H^{n+1}(N_i,N_{i-1})\to 0\\
 &\downarrow     &                  &\searrow &              \\
 &  0                   &                  &         & H^{n+1}(N_i,N_i').
  \end{array}
\end{equation*}
which implies 
%\begin{subequations}\label{top}
\begin{equation}
H^n(N_i,N_{i-1})=\ker\Delta_i={\cB}_i(\coker{\cD}_i\cap\ker{\cD}_i^*)
\end{equation}
\begin{equation}
H^{n+1}(N_i,N_{i-1})=\coker\Delta_i=\coker{\cD}_i^*= H^{n+1}(N_i,N_i')
\end{equation}
%\end{subequations}
Sequence \eqref{minima} has torsion $\tau({\cD}_i)^{(-1)^{n-1}}$
and sequence \eqref{maxima} has torsion $\tau(\Delta_i)^{(-1)^{n-1}}$. 
Therefore, Corollary \ref{aditive} gives
%\begin{subequations}\label{torsion-odd}
\begin{equation}
\tau_{d(N_i'',N_{i-1})}=\tau_{d(N_i'',N_i')}\tau_{d(N_i',N_{i-1})}
\tau({\cD}_i)^{(-1)^n}=\tau({\cD}_i)^{(-1)^n},
\end{equation}
\begin{equation}\label{torsion-odd}
\tau_{d_0^{n-1,i}}=\tau_{d(N_i,N_i'')}\tau_{d(N_i'',N_{i-1})}
\tau(\Delta_i)^{(-1)^{n-1}}=\tau({\cD}_i)^{(-1)^n}\tau(\Delta_i)^{(-1)^{n-1}}.
\end{equation}
%\end{subequations}
To compute $\tau({\cD}_i)$, we recall that 
${\cD}_i=\left(\begin{smallmatrix}
I-\rho(\alpha_i)\\I\pm\rho(\beta_i)\end{smallmatrix}\right)$.
Since $\rho(\alpha_i), \rho(\beta_i)\in\U(\m)$ and they commute. There is
a splitting 
$${\C}^m=V_{\alpha_i}+V_{\beta_i}+\ker{\cD}_i,$$
invariant under $\rho(\alpha_i)$, and $\rho(\beta_i)$, such that
$I-\rho(\alpha_i):V_{\alpha_i}\hookleftarrow$ and 
$I-\rho(\beta_i):V_{\beta_i}\hookleftarrow$ are nonsingular. Let 
${\bf k}_i$, ${\bf v}_{\alpha_i}$, ${\bf v}_{\beta_i}$ be bases of 
$H^{n-1}(N_i,N_{i-1})=\ker{\cD}_i$, $V_{\alpha_i}$, $V_{\beta_i}$. Let
$$W_i=(V_{\beta_i}+\ker{\cD}_i)\oplus (V_{\alpha_i}+\ker{\cD}_i),$$
then $j:W_i\to H^n(N_i'',N_{i-1})$ is an isomorphism, 
and so $j[({\bf v}_{\beta_i}\cup {\bf k}_i)\times\{0\}]
\cup j[\{0\}\times({\bf v}_{\alpha_i}\cup {\bf k}_i)]$ is a basis of 
$H^n(N_i'',N_{i-1})$. Thus
\begin{equation}\label{torsion-Di}
\begin{split}
\tau({\cD}_i) & =
\frac{\wedge({\cD}_i{\bf v}_{\alpha_i}, {\cD}_i{\bf v}_{\beta_i},
({\bf v}_{\beta_i}\cup {\bf k}_i)\times\{0\},\{0\}
\times({\bf v}_{\alpha_i}\cup {\bf k}_i))\otimes\wedge{\bf k}_i}
{\wedge({\bf v}_{\alpha_i},{\bf v}_{\beta_i},{\bf k}_i)
\otimes\wedge ({\cB}_i[({\bf v}_{\beta_i}\cup {\bf k}_i)\times\{0\}],
{\cB}_i[\{0\}\times({\bf v}_{\alpha_i}\cup {\bf k}_i)])}\\
= & \frac{\wedge((I-\rho(\alpha_i)){\bf v}_{\alpha_i},{\bf v}_{\beta_i},
{\bf k}_i)\otimes\wedge({\bf v}_{\alpha_i},(I-\rho(\beta_i)){\bf v}_{\beta_i},
{\bf k}_i)\otimes\wedge{\bf k}_i}
{\wedge({\bf v}_{\alpha_i},{\bf v}_{\beta_i},{\bf k}_i)\otimes 
{\cB}_i^*(\wedge ({\bf v}_{\beta_i},{\bf k}_i)\otimes
\wedge({\bf v}_{\alpha_i},{\bf k}_i))}.
\end{split}
\end{equation}
To compute $\tau(\Delta_i)$ we recall that 
${\cD}_i^*=(I\pm\rho(\beta_i),\rho(\alpha_i)-I)$ and then
$$W_i\cap \ker{\cD}_i^*=\ker{\cD}_i\oplus\ker{\cD}_i,\,
{\cD}_i^*({\C}^m\oplus{\C}^m)={\cD}_i^*(W_i)=V_{\beta_i}+V_{\alpha_i}.$$
Thus $\ker\Delta_i={\cB}_i(\ker{\cD}_i\oplus\ker{\cD}_i)$ and
$\coker\Delta_i=\coker{\cD}_i^*=\ker{\cD}_i$. 
Thus  $H^n(N_i,N_{i-1})\cong\ker{\cD}_i\oplus\ker{\cD}_i$,
$H^n(N_i,N_{i-1})=\ker{\cD}_i$, and using  ${\bf k}_i\times{\bf k}_i$
and ${\bf k}_i$ as their bases we have
\begin{equation}\label{torsion-delta}
\tau(\Delta_i)=\frac{\wedge((I-\rho(\alpha_i)){\bf v}_{\alpha_i},
(I-\rho(\beta_i)){\bf v}_{\beta_i},{\bf k}_i)
\otimes(\wedge{\bf k}_i\otimes\wedge{\bf k}_i)}
{{\cB}_i^*(\wedge ({\bf v}_{\beta_i},{\bf k}_i)\otimes 
\wedge({\bf v}_{\alpha_i},{\bf k}_i))\otimes\wedge{\bf k}_i}.
\end{equation}
From equations \eqref{torsion-odd}, \eqref{torsion-Di},
\eqref{torsion-delta} we have
$$\tau_{d_0^{n-1,i}}=1.$$
\end{pf}
We now come to the description of the components $F_{ij}^*$ of $d_1$. 
Let $\psi_t$ be the gradient flow of $g$.
One can construct an index filtration 
\begin{equation*}
\emptyset=K_{-1}\subset K_0\subset L_1\subset P_1\subset K_1\subset L_2
\subset P_2\subset K_2\subset K_3=M
\end{equation*}
such that $L_1\subset M_1\subset L_2$, 
$P_1\subset M_2\subset P_2$ and
\begin{equation*}
\begin{split}
\bigcap_{t\in{\R}}\psi_t(K_i\setminus K_{i-1})&=C_i(g)\quad i=0,1,2,3.\\
\bigcap_{t\in{\R}}\psi_t(L_i\setminus K_{i-1})&=C_i(g)\cap M_1\quad i=1,2.\\
\bigcap_{t\in{\R}}\psi_t(P_i\setminus L_i)&=C_i(g)\cap(M_2\setminus M_1)
\quad i=1,2.\\
\bigcap_{t\in{\R}}\psi_t(K_i\setminus P_i)&=C_i(g)\setminus M_2\quad i=1,2.
\end{split}
\end{equation*}
We have
\begin{equation}\label{cohom-K}
  H^i(K_i,K_{i-1})=\bigoplus_{\text{index}(x)=i}{\C}^mx.
\end{equation}
Note that all the components $G^*_{xy}:{\C}^mx\to{\C}^my$
of the maps $d^K$ in the cochain complex
\begin{equation}\label{cochain-K}
0\to  H^0(K_0)\buildrel d^K\over\longrightarrow H^1(K_1,K_0)\buildrel d^K\over
\longrightarrow H^2(K_2,K_1)\buildrel d^K\over\longrightarrow H^2(M,K_2)\to 0,
\end{equation}
are given as in equation \eqref{Floer} of Lemma\ref{Witten}.

\begin{theorem}\label{first-diff}
For critical circles $\gamma_i$ the components $F^k_{ij}$ are induced by the 
maps 
\begin{itemize}
\item(c.c) $G^n_{w_iw_j}$, $G^{n+1}_{z_iz_j}$ for $k_{2n}<i\le k_{2n+1}$, 
$k_{2n+2}<j\le k_{2n+3}$, $n=0,1$
\item(c.t) $G^1_{w_iq_j}, G^1_{w_ir_j}, G^2_{z_is_j}$ 
for $k_2<i\le k_3<j\le k_4$.
\end{itemize}
For critical tori $F_{c_i}$ the components $F^k_{ij}$ are induced
by the maps
\begin{itemize}
\item(t.c) $G^0_{p_iw_j}$, $G^1_{q_iz_j}$, $G^1_{r_iz_j}$ for 
$k_1<i\le k_2<j\le k_3$.
\end{itemize}
\end{theorem}
\begin{pf}
Consider the commutative diagrams
\begin{equation}\label{dM1M2}
  \begin{array}{ccccc}
H^0(M_1)& \buildrel d_1\over\longrightarrow & H^1(M_2,M_1) &
\buildrel d_1\over\longrightarrow & H^2(M,M_2)\\
& & \downarrow     &  &\uparrow  \\
& &H^1(M_2,L_1)&  & H^2(M,P_2)\\ 
\downarrow & &\downarrow & &\downarrow \\     
& & H^1(P_1,L_1)&  & H^2(K_2,P_2)\\
& &\downarrow  & & \\
H^0(K_0)& & H^1(P_1,K_0) & &\downarrow \\
&\searrow d^K & \uparrow & & \\
& & H^1(K_1,K_0)& \buildrel d^K\over\longrightarrow & H^2(K_2,K_1)
  \end{array}
\end{equation} 
and
\begin{equation}\label{dM2M3}
  \begin{array}{ccccc}
H^1(M_1)&\buildrel d_1\over\longrightarrow & H^2(M_2,M_1) 
&\buildrel d_1\over\longrightarrow & H^3(M,M_2)\\
 \downarrow & &\uparrow & &\\
H^1(L_1) & & H^2(P_2,M_1)& & \\
\uparrow & &\uparrow & &\uparrow\\
H^1(L_1,K_0)& &H^2(P_2,K_1)& & \\
\downarrow & &\downarrow & &\\
H^1(K_1,K_0)& & H^2(P_2,L_2)& &H^3(M,K_2)\\
& \searrow d^K &\downarrow &d^K \nearrow &\\
& & H^2(K_2,K_1)& & 
  \end{array}
\end{equation}
where the downwards maps are injective and the upwards maps are 
surjective. By \eqref{dM1M2} the maps $d_1$ in
\begin{equation*}
H^0(M_1)\buildrel d_1\over\longrightarrow H^1(M_2,M_1)
\buildrel d_1\over\longrightarrow H^2(M,M_2)
\end{equation*}
are induced by the components  $H^0(K_0)\to H^1(P_1,L_1)$ and 
$$H^1(P_1,L_1)=\bigoplus_{k_2<i\le k_3}{\C}^mw_i\to H^2(K_2,P_2)=
\bigoplus_{k_4<i}{\C}^mw_i\oplus
\bigoplus_{k_3<i\le k_4}{\C}^mq_i\oplus{\C}^mr_i$$
of the maps $d^K$ in \eqref{cochain-K}. By \eqref{dM2M3} the maps $d_1$ in
\begin{equation*}
H^1(M_1)\buildrel d_1\over\longrightarrow H^2(M_2,M_1)
\buildrel d_1\over\longrightarrow H^3(M,M_2)
\end{equation*}
are induced by the components 
$$H^1(L_1,K_0)=
\bigoplus_{i\le k_1}{\C}^mz_i\oplus
\bigoplus_{k_1<i\le k_2}{\C}^mq_i\oplus{\C}^mr_i\to H^2(P_2,L_2)=
\bigoplus_{k_2<i\le k_3}{\C}^mz_i$$
and $H^2(P_2,L_2)\to H^3(M,K_2)$
of the maps $d^K$ in \eqref{cochain-K}.
\end{pf}
We now consider the term $(E_2,d_2)$ of the spectral sequence defined
by the filtration. By \eqref{E2} and Corollary \ref{E1}, the spaces $E_2^{n,q}$
that can be nonzero are
$$E_2^{0,q}=\ker(d_1:H^q(M_1)\to
H^{q+1}(M_2,M_1)),\,q=0,1.E_2^{0,2}=H^2(M_1)$$
$$E_2^{1,q}=\ker(d_1:H^{q+1}(M_2,M_1)\to H^{q+2}(M,M_2))/d_1(H^q(M_1)),\,
q=0,1$$
$$E_2^{2,-1}=H^1(M,M_2),\,E_2^{2,q}=H^{q+2}(M,M_2)/d_1(H^{q+1}(M_2,M_1)),\,
q=0,1.$$

Therefore, we can have nonzero maps $d_2:E_2^{n,q}\to E_2^{n+2,q-1}$ only
for $n=0$, $q=0,1,2$.

\begin{theorem}\label{second-diff}\quad

\begin{itemize}
\item[(0)] The map $d_2:E_2^{0,0}\to E_2^{2,-1}$ is induced by the component
\newline
$H^0(K_0)\to\underset{k_3<i\le k_4}{\bigoplus}{\C}^mp_i$
of $d^K:H^0(K_0)\to H^1(K_1,K_0)$.
\item[(1)] The map $d_2:E_2^{0,1}\to E_2^{2,0}$ is induced by the component
$$\bigoplus_{i\le k_1}{\C}^mz_i\oplus\bigoplus_{k_1<i\le k_2}
{\C}^mq_i\oplus{\C}^mr_i\to\bigoplus_{k_4<i}{\C}^mw_i\oplus
\bigoplus_{k_3<i\le k_4}{\C}^mq_i\oplus{\C}^mr_i$$
of $d^K:H^1(K_1,K_0)\to H^2(K_2,K_1)$.
\item[(2)] The map $d_2:E_2^{0,2}\to E_2^{2,1}$ is induced by the component 
\newline
$\underset{k_1<i\le k_2}{\bigoplus}{\C}^ms_i\to H^3(M,K_2)$ of
$d^K:H^2(K_2,K_1)\to H^3(M,K_2)$.
\end{itemize}
\end{theorem}
\begin{pf}
By \eqref{eq:second-diff}, the map $d_2:E_2^{0,q}\to E_2^{2,q-1}$ is
given as the composition 
$${\j}_q^{-1}\circ\delta_q|\ima k_q:
\ima k_q\to H^{q+1}(M,M_2)/\ima j_q$$
where $k_q:H^q(M_2)\to H^(M_1)$, $\delta:H^q(M_1)\to H^{q+1}(M,M_1)$ is the 
coboundary map and $j_q:H^{q+1}(M,M_2)\to H^{q+1}(M,M_1)$ defines the
isomorphism ${\j}_q:H^{q+1}(M,M_2)/\ker j_q\to\ima j_q$.
Consider the commutative diagrams
\begin{equation}\label{d:E20}
  \begin{array}{ccccc}
E_2^{0,0} & \buildrel\delta_0|\ima k_0\over\longrightarrow & \ima j_0
& \buildrel{\j}_0\over\longleftarrow & E_2^{2,-1}\\
\cap & & \cap & & \uparrow\\
H^0(M_1) & \buildrel\delta_0\over\longrightarrow & H^1(M,M_1) &
\buildrel j_0\over\longleftarrow & H^1(M,M_2)\\
& & & &\downarrow\\
& & & &H^1(M,P_1)\\
\downarrow & & & & \downarrow\\
& & & &H^1(K_1,P_1)\\
& & & &\downarrow\\
H^0(K_0) & &\buildrel d^K\over\longrightarrow  & & H^1(K_1,K_0)
\end{array}
\end{equation}
\begin{equation}\label{d:E21}
  \begin{array}{ccccc}
E_2^{0,1} & \buildrel\delta_1|\ima k_1\over\longrightarrow & \ima j_1 
& \buildrel{\j}_1\over\longleftarrow & E_2^{2,0}\\ 
\cap & & \cap & & \uparrow\\
H^1(M_1)&\buildrel\delta_1\over\longrightarrow & H^2(M,M_1) & 
\buildrel j_1\over\longleftarrow &  H^2(M,M_2)\\
\downarrow & & & & \uparrow d^1\\
H^1(L_1) & & & & H^2(M,P_2)\\
\uparrow & & & &\downarrow\\
H^1(L_1,K_0)& & & &H^2(K_2,P_2)\\
\downarrow & & & &\downarrow \\
H^1(K_1,K_0)& & \buildrel d^K\over\longrightarrow & & H^2(K_2,K_1)
  \end{array}
\end{equation}
\begin{equation}\label{d:E22}
  \begin{array}{ccccc}
E_2^{0,2} & \buildrel\delta_2|\ima k_2\over\longrightarrow & \ima j_2 & 
\buildrel{\j}_0\over\longleftarrow & E_2^{2,1}\\
\cap & & \cap & & \uparrow\\
H^2(M_1)& \buildrel\delta\over\longrightarrow & H^3(M,M_1) & 
\buildrel j\over\longleftarrow &  H^3(M,M_2)\\
\uparrow & & & &\\
H^2(L_2) & & & & \\
\uparrow & & & &\uparrow \\
H^2(L_2,K_1)& && & \\
\uparrow & & & & \\
H^2(K_2,K_1) & & \buildrel d^K\over\longrightarrow & & H^3(M,K_2)
\end{array}
\end{equation}
where as in Theorem \ref{first-diff} the downwards maps are injective
and the upwards maps are surjective. By \eqref{d:E20} and 
$H^1(K_1,P_1)=\underset{k_3<i\le k_4}{\bigoplus}{\C}^mp_i$, we get (0).
By \eqref{d:E21}, 
$H^1(L_1,K_0) =\underset{i\le k_1}{\bigoplus}{\C}^mz_i\oplus
\underset{k_1<i\le k_2}{\bigoplus}{\C}^mq_i\oplus{\C}^mr_i$ and
$H^2(K_2,P_2)=\underset{k_4<i}{\bigoplus}{\C}^mw_i\oplus
\underset{k_3<i\le k_4}{\bigoplus}{\C}^mq_i\oplus{\C}^mr_i$ we get (1). By 
\eqref{d:E21} and $H^2(L_2,K_1)=\underset{k_1<i\le k_2}{\bigoplus}{\C}^ms_i$, 
we get (2).
\end{pf}

\section{Examples}
Consider the following instance of the Kowalevskaya integrable case
of the rigid body.
There are two minimal circles $m_1, m_2$, two nonorientable and one 
orientable saddle circles $r_1, r_2, r_3$, and one
maximal circle $n$. The family of tori starting at $m_i$ changes
to a family ${\cF}_i$ of tori when crossing $r_i$ $(i=1,2)$. The 
families ${\cF}_1$ and ${\cF}_2$ come together to become one family
when crossing $r_3$. The manifold $M$ is the 3-dimensional real 
projective space and so $\pi_1(M)={\Z}_2$. A representation 
$\rho:\pi_1(M)\to{\U}(1)$ is given by $\rho([0])=1$, $\rho([1])=-1$.
%\centerline{\epsfysize=3cm\epsfbox{ochos.eps}}
We have that $\rho(m_i)=1$, $\rho(r_j)=\rho(n)=-1$,
$\Delta(r_i)=-1$, $\Delta(r_3)=1$, $i=1,2,\,j=1,2,3$. Therefore
%\begin{subequations}
\begin{equation}
H^k(M_1)=\begin{cases} {\C}\oplus{\C}&\text{if}\quad k=0,1\\
0&\text{in other case}\end{cases}
\end{equation}
\begin{equation}
H^k(M_2,M_1)=\begin{cases} {\C}\oplus{\C}&\text{if}\quad k=1,2\\
0&\text{in other case}\end{cases}
\end{equation}
\begin{equation} 
H^*(M,M_2)=0
\end{equation}
\begin{equation}
\tau_{d_0^{0,i}}=\tau_{d_0^{1,i}}=1\,(i=1,2),\,\tau_{d_0^{1,3}}=\frac 12,
\tau_{d_0^2}=2\Rightarrow \tau_{d_0}=1
\end{equation}
%\end{subequations}

We now change the Bott integral $f$ to a Morse Smale function $g$ 
to be able to compute $d_1:H^k(M_1)\to H^{k+1}(M_2,M_1)$, $k=0,1$.
Note that $W^u(r_i),i=1,2$ is a M\"obius strip $\Sigma_i$ with
$\partial\Sigma_i=m_i$. The function $g$ has critical points $w_i,z_i$
on $m_i$ with indices $0,1$ respectively, and critical points
$\eta_i,\zeta_i$ on $r_i$ with indices $1,2$ respectively. There
is one orbit $\alpha_i$ connecting $\zeta_i$ to $z_i$ and two orbits
$\beta_i,\delta_i$ connecting $\eta_i$ to $w_i$, with 
$\rho(\alpha_i)=1$, $\rho(\beta_i)=1$ $\rho(\delta_i)=-1$. Thus
$G_{w_i\eta_i}:{\C}\to{\C}=2$ and 
$G_{z_i\zeta_i}:{\C}\to{\C}=1$. 

Therefore $d_1:H^0(M_1)\to H^1(M_2,M_1)=
\bigl(\begin{smallmatrix}2&0\\0&2\end{smallmatrix}\bigr)$ and
$d_1:H^1(M_1)\to H^2(M_2,M_1)=
\bigl(\begin{smallmatrix}1&0\\0&1\end{smallmatrix}\bigr)$. Thus
$E_2=0$ and $\tau_{d_1}=2\cdot 2/1=4$.

Institut f\"ur  Mathematik,E.-M.-Arndt- Universit\"at  Greifswald 
 
Jahn-strasse 15a, D-17489 Greifswald, Germany.
 
{\it E-mail address}: felshtyn@@rz.uni-greifswald.de

Instituto de Mathematicas, Universidad Nacional Autonoma de Mexico,
Ciudad Universitaria C.P. 04510, Mexico D.F., Mexico.

{\it E-mail address}: hector@@math.unam.mx

\end{document}